\documentclass[12pt]{article}
\usepackage{amsfonts, amssymb}

\def\AA{{\bf A}}
\def\E{{\mathcal E}}
\def\F{{\mathcal F}}
\def\G{{\mathcal G}}
\def\FF{{\bf F}}

\def\mm{{\mathfrak m}}
\def\O{{\mathcal O}}
\def\PP{{\bf P}}
\def\Q{\mathbb Q}
\def\Qq{\mathfrak Q}
\def\Z{\mathbb Z}
\def\dim{\mathop{\rm dim}\nolimits}

\def\spec{\mathop{\rm Spec}\nolimits}
\def\sym{\mathop{\rm Sym}\nolimits}

\def\hom{{\mathfrak Hom}}

\let\hra\hookrightarrow

\newtheorem{theorem}{Theorem}
\newtheorem{corollary}[theorem]{Corollary}
\newtheorem{lemma}[theorem]{Lemma}
\def\proof{\paragraph{Proof}}
\makeatletter \def\l@section{\@dottedtocline{1}{0em}{1.2em}} \makeatother

\begin{document}

\title{Representability of Hom implies flatness}

\author{Nitin Nitsure}
 
\date{}

\maketitle

\begin{abstract}
Let $X$ be a projective scheme over a noetherian base scheme $S$, and 
let $\F$ be a coherent sheaf on $X$. For any coherent sheaf $\E$ on $X$, 
consider the set-valued contravariant functor $\hom_{(\E,\F)}$ on 
$S$-schemes, defined by $\hom_{(\E,\F)}(T) = Hom(\E_T,\F_T)$ where 
$\E_T$ and $\F_T$ are the pull-backs of $\E$ and $\F$ to 
$X_T = X\times_S T$. A basic result of Grothendieck ([EGA] III 7.7.8, 7.7.9) 
says that if $\F$ is flat over $S$ then $\hom_{(\E,\F)}$ is representable 
for all $\E$. 

We prove the converse of the above, in fact, we show that if $L$ 
is a relatively ample line bundle on $X$ over $S$ such that the 
functor $\hom_{(L^{-n},\F)}$ is representable for infinitely many 
positive integers $n$, then $\F$ is flat over $S$. 
As a corollary, taking $X=S$, it follows that if $\F$ is a coherent sheaf 
on $S$ then the functor $T\mapsto H^0(T, \F_T)$ on the category of 
$S$-schemes is representable if and only if $\F$ is locally free on $S$. 
This answers a question posed by Angelo Vistoli. 

The techniques we use 
involve the proof of flattening stratification, together with the methods 
used in proving the author's earlier result (see [N1]) that the 
automorphism group functor of a coherent sheaf on $S$ is representable 
if and only if the sheaf is locally free. 
\end{abstract}

\bigskip

Let $S$ be a noetherian scheme, and let $X$ be a projective scheme
over $S$. If $\E$ and $\F$ are coherent sheaves on $X$, consider the 
contravariant functor $\hom_{(\E,\F)}$ from the category of 
schemes over $S$ to the category of sets
which is defined by putting
$$\hom_{(\E,\F)}(T) = Hom_{X_T}(\E_T,\F_T)$$
for any $S$-scheme $T\to S$, where $X_T = X\times_ST$, and 
$\E_T$ and $\F_T$ denote the pull-backs of $\E$ and $\F$ under the
projection $X_T \to X$. 
This functor is clearly a sheaf
in the fpqc topology on $Sch/S$.
It was proved by Grothendieck that if $\F$ is flat over
$S$ then the above functor is representable (see [EGA] III 7.7.8, 7.7.9).

Our main theorem is as follows, which is a converse to the above. 

\begin{theorem}\label{main theorem} 
Let $S$ be a noetherian scheme, $X$ a projective scheme
over $S$, and $L$ a relatively very ample line bundle on $X$ over $S$.
Let $\F$ be a coherent sheaf on $X$. Then the following three 
statements are equivalent:

(1) The sheaf $\F$ is flat over $S$.

(2) For any coherent sheaf $\E$ on $X$, the set-valued contravariant functor 
$\hom_{(\E,\F)}$ on $S$-schemes, defined by 
$\hom_{(\E,\F)}(T) = Hom_{X_T}(\E_T,\F_T)$, is representable.

(3) There exist infinitely many positive integers $r$ such that 
the set-valued contravariant functor $\G^{(r)}$ on $S$-schemes, 
defined by $\G^{(r)}(T) =  H^0(X_T,\F_T\otimes L^{\otimes r})$, is 
representable.
\end{theorem}

In particular, taking $X=S$ and $L = \O_X$, 
we get the following corollary.

\begin{corollary}\label{Question of Vistoli}
Let $S$ be a noetherian scheme, and $\F$ a coherent sheaf on $S$.
Consider the contravariant functor $\FF$ from $S$-schemes to sets,
which is defined by putting $\FF(T) = H^0(T, f^*\F)$ for any 
$S$-scheme $f: T\to S$. This functor (which is 
a sheaf in the fpqc topology) is representable if and only if
$\F$ is locally free as an $\O_S$-module.
\end{corollary}

Note that the affine line $\AA^1_S$ over a base $S$ admits a  
ring-scheme structure over $S$ in the obvious way. 
A {\bf linear scheme} over a scheme $S$ will mean a module-scheme
$V\to S$ under the ring-scheme $\AA^1_S$. This means $V$ is 
a commutative group-scheme over $S$ together with 
a `scalar-multiplication' morphism
$\mu : \AA^1_S \times_S V \to V$ over $S$, such that 
the module axioms (in diagrammatic terms) are satisfied.

A {\bf linear functor} $\FF$ on $S$-schemes will mean a 
contravariant functor from $S$-schemes to sets 
together with the structure of an $H^0(T,\O_T)$-module 
on $\FF(T)$ for each $S$-scheme $T$, 
which is well-behaved under any
morphism $f: U\to T$ of $S$-schemes in the following sense: 
$\FF(f): \FF(T) \to \FF(U)$ is a homomorphism of the 
underlying additive groups, and 
$\FF(f)(a\cdot v) = f^*(a)\cdot (\FF(f) v)$
for any
$a\in H^0(T,\O_T)$ and $v \in \FF(T)$. 
In particular note that the kernel of $\FF(f)$ 
will be an $H^0(T,\O_T)$-submodule of $\FF(T)$. 
The functor of points of a linear scheme is naturally a linear functor.
Conversely, it follows by the Yoneda lemma
that if a linear functor $\FF$ on $S$-schemes is representable, then
the representing scheme $V$ is naturally a linear scheme over $S$.

For example, the linear functor $T\mapsto H^0(T,\O_T)^n$ (where $n\ge 0$)
is represented by the affine space $\AA^n_{\Z}$ over $\spec \Z$,
with its usual linear-scheme structure. More generally, for any
coherent sheaf $\Qq$ on $S$, the scheme ${\bf Spec} \sym(\Qq)$
is naturally a linear-scheme over $S$, where $\sym(\Qq)$ denotes the 
symmetric algebra of $\Qq$ over $\O_S$. It represents the 
linear functor $\FF(T) =  Hom(\Qq_T,\O_T)$ where $\Qq_T$ 
denotes the pull-back of $\Qq$ under $T\to S$.

With this terminology, the functor 
$\G^{(r)}(T) =  H^0(X_T,\F_T\otimes L^{\otimes r})$ of 
Theorem \ref{main theorem}.(3) is a linear functor. Therefore, if 
a representing scheme $G^{(r)}$ exists, 
it will naturally be a linear scheme. Note that each $\G^{(r)}$ is
obviously a sheaf in the fpqc topology.

The proof of Theorem \ref{main theorem} is by a combination of 
the result of Grothendieck on the existence of a 
flattening stratification ([TDTE -IV]) together
with the techniques which were employed in [N1] 
to prove the following result.

\begin{theorem}\label{GL}
{\bf (Representability of the functor $GL_E$) }  
Let $S$ be a noetherian scheme, and $E$ a coherent $\O_S$-module. 
Let $GL_E$ denote the contrafunctor on $S$-schemes
which associates to any $S$-scheme $f:T\to S$ 
the group of all $\O_T$-linear
automorphisms of the pull-back $E_T= f^*E$ 
(this functor is a sheaf in the fpqc topology).
Then $GL_E$ is representable by a group scheme over $S$ 
if and only if $E$ is locally free.
\end{theorem}

We re-state Grothendieck's result (see [TDTE IV]) on 
the existence of a flattening stratification in the following
form, which emphasises the role of the direct images $\pi_*(\F(r))$. 
For an exposition of flattening stratification, see Mumford [M] or [N2]. 

\begin{theorem}\label{flattening stratification}{\rm (Grothendieck)} 
Let $S$ be a noetherian scheme, and let $\F$ be a coherent sheaf on $\PP^n_S$
where $n\ge 0$. There exists an integer $m$, and 
a collection of locally closed subschemes $S_f \subset S$ indexed
by polynomials $f \in \Q[\lambda]$, with the following properties.

(i) The underlying set of $S_f$ consists of all $s\in S$ such that 
the Hilbert polynomial of $\F_s$ is $f$, where $\F_s$ denotes the 
pull-back of $\F$ to the schematic fibre $\PP^n_s$ over $s$ of the 
projection $\pi : \PP^n_S \to S$. All but finitely many
$S_f$ are empty (only finitely many Hilbert polynomials occur). 
In particular, the $S_f$ are mutually disjoint, and their set-theoretic union 
is $S$.

(ii) For each $r\ge m$, the higher direct images $R^j\pi_*(\F(r))$
are zero for $j\ge 1$ and 
the subschemes $S_f$ give the 
flattening stratification for the direct image
$\pi_*(\F(r))$, that is, the morphism $i: \coprod_f S_f \to S$
induced by the locally closed embeddings 
$S_f \hra S$ has the universal property that 
for any morphism $g : T\to S$, the sheaf $g^*\pi_*(\F(r))$ is locally free
on $T$ if and only if $g$ factors via $i: \coprod_f S_f \to S$.

(iii) The subschemes $S_f$ give the 
flattening stratification for $\F$, that is,
for any morphism $g : T\to S$, the sheaf $\F_T = (1\times g)^*\F$
on $\PP^n_T$ is flat over $T$ if and only if $g$ 
factors via $i: \coprod_f S_f \to S$. In particular, 
$\F$ is flat over $S$ if and only if each $S_f$ is an open
subscheme of $S$. 

(iv) Let $\Q[\lambda]$ be totally ordered by putting $f_1 < f_2$ 
if $f_1(p) < f_2(p)$ for all $p \gg 0$. Then the closure of 
$S_f$ in $S$ is set-theoretically contained in $\bigcup_{g\ge f}\,S_g$.
Moreover, whenever $S_f$ and $S_g$ are non-empty, we have
$f<g$ if and only if $f(p) < g(p)$ for all $p\ge m$.
\end{theorem}

The following elementary lemma of Grothendieck on base-change  
does not need any flatness hypothesis.
The price paid is that the integer $r_0$ may depend on $\phi$.
(See [N2] for a cohomological proof.)

\begin{lemma}\label{base change without flatness} 
Let $\phi : T\to S$ be a morphism of noetherian schemes, 
let $\F$ a coherent sheaf on $\PP^n_S$, and 
let $\F_T$ denote its pull-back 
under the induced morphism $\PP^n_T\to \PP^n_S$. Let
$\pi_S : \PP^n_S \to S$ and $\pi_T : \PP^n_T \to T$
denote the projections. Then there exists an integer 
$r_0$ such that the base-change homomorphism 
$\phi^* {\pi_S}_*\, \F(r) \to {\pi_T}_* \, \F_T(r)$
is an isomorphism for all $r\ge r_0$.
\end{lemma}

\bigskip

{\bf Proof of Theorem \ref{main theorem} } 
The implication $(1) \Rightarrow (2)$ follows by [EGA] III 7.7.8, 7.7.9,
while the implication $(2) \Rightarrow (3)$ follows by taking 
$\E = L^{\otimes -r}$. Therefore it now remains to show 
the implication $(3) \Rightarrow (1)$. This we do in a number of steps.

\medskip

{\bf Step 1: Reduction to $S=\spec R$ with $R$ local, $X=\PP^n_S$ and 
$L= \O_{\PP^n_S}(1)$}\\
Suppose that $\F$ is not flat over $S$, but 
the linear functor $\G^{(r)}$ on $S$-schemes, defined by
$\G^{(r)}(T) =  H^0(X_T,\F_T\otimes L^{\otimes r})$, is 
representable by a linear scheme $G^{(r)}$
over $S$ for arbitrarily large integers $r$. 
As $\F$ is not flat, by definition there exists 
some $x\in X$ such that the stalk $\F_x$ is not a flat module over
the local ring $\O_{S,\pi(x)}$ where 
$\pi: X\to S$ is the projection. Let $U = \spec \O_{S,\pi(x)}$,
let $\F_U$ be the pull-back of $\F$ to $X_U = X\times_SU$
and let $G^{(r)}_U$ denote the pull-back of $G^{(r)}$ to $U$. 
Then $\F_U$ is not flat over $U$ but 
given any integer $m$, there exists an integer $r\ge m$ such that 
the functor $\G^{(r)}_U$ on $U$-schemes, defined by
$\G^{(r)}_U(T) =  H^0(X_T,\F_T\otimes L^{\otimes r})$, is 
representable by the $U$-scheme $G^{(r)}_U$. 

Therefore, by replacing $S$ by $U$, we can assume that $S$ is of the 
form $\spec R$ where $R$ is a noetherian local ring. 
Let $i: X\hra \PP^n_S$ be the embedding given by $L$. 
Then replacing $\F$ by $i_*\F$, 
we can further assume that $X=\PP^n_S$ and $L= \O_{\PP^n_S}(1)$.

\medskip

{\bf Step 2: Flattening stratification of $\spec R$ }
There exists an integer $m$ as asserted by
Theorem \ref{flattening stratification}, such that for any $r\ge m$,
the flattening
stratification of $S$ for the sheaf $\pi_*\F(r)$ on $S$ is the 
same as the flattening
stratification of $S$ for the sheaf $\F$ on $\PP^n_S$.
Let $r\ge m$ be any integer.
As $\F$ is not flat over $S = \spec R$, 
the sheaf $\pi_*\F(r)$ is not flat.
Let $M_r = H^0(S, \pi_*\F(r))$, which is a finite $R$-module.
Let $\mm\subset R$ be the maximal ideal, and let $k = R/\mm$ the residue field.
Let $s\in S = \spec R$ be the closed point, and let 
$d = \dim_k(M_r/\mm M_r)$. Then there exists a right-exact sequence
of $R$-modules of the form
$$R^{\delta} \stackrel{\psi}{\to} R^d \to M_r\to 0$$
Let $I\subset R$ be the ideal formed by the matrix entries of 
the $(d\times \delta)$-matrix $\psi$. Then $I$ defines a
closed subscheme $S'\subset S$ which is the flattening stratification of
$S$ for $M_r$. As $M_r$ is not flat by assumption, $I$ is a non-zero
proper ideal in $R$.

It follows from Theorem \ref{flattening stratification} that   
$I$ is independent of $r$ as long as $r\ge m$.

\medskip

{\bf Step 3: Reduction to artin local case with principal $I$ with $\mm I =0$ }
Let $I = (a_1,\ldots, a_t)$ where $a_1,\ldots, a_t$ is a 
minimal set of generators of $I$. Let 
$J \subset R$ be the ideal defined by
$$J = (a_2,\ldots,a_t) + \mm I$$
Then note that $J \subset I\subset \mm$, and the quotient
$R' = R/J$ is an artin local $R$-algebra with maximal ideal
$\mm' = \mm/J$, and $I' = I/J$ is a non-zero
principal ideal which satisfies $\mm' I' =0$.   
For the base-change under $f: \spec R' \to \spec R$, 
the flattening stratification for $f^*\pi_*\F(r)$ is defined by the 
ideal $I'$ for $r\ge m$. 
Let $\F'$ denote the pull-back of $\F$ to $\PP^n_{R'}$, and let
$\pi' : \PP^n_{R'} \to \spec R'$ the projection.
As $f$ is a morphism of noetherian schemes, by 
Lemma \ref{base change without flatness} 
there exists some integer $m'$ such that the base-change
homomorphism
$f^*\pi_*\F(r) \to \pi'_*\F'(r)$
is an isomorphism whenever $r\ge m'$.
Choosing some $m' \ge m$ with this property, and replacing
$R$ by $R'$, $\F$ by $\F'$ and $m$ by $m'$, 
we can assume that $R$ is artin local, 
and $I$ is a non-zero principal ideal with $\mm I =0$,
which defines the flattening stratification for $\pi_*\F(r)$
for all $r\ge m$.

\medskip

\newpage

{\bf Step 4: Decomposition of $\pi_*\F(r)$ via lemma of Srinivas }

{\bf Lemma }(Srinivas)
{\it Let $R$ be an artin local ring with maximal ideal $\mm$, and let
$E$ be any finite $R$ module whose flattening stratification 
is defined by an ideal $I$ which is a 
non-zero proper principal ideal with $\mm I =0$.
Then there exist integers $i\ge 0$ and $j>0$ such that 
$E$ is isomorphic to the direct sum $R^i\oplus (R/I)^j$.
}

{\bf Proof } See Lemma 4 in [N1].

\medskip

We apply the above lemma to the 
$R$-module $M_r = H^0(S, \pi_*\F(r))$, which has flattening
stratification defined by the principal ideal $I$ with $\mm I =0$,
to conclude that (up to isomorphism) $M_r$ has the form
$$M_r = R^{i(r)}\oplus (R/I)^{j(r)}$$
for non-negative integers $i(r)$ and $j(r)$ with $j(r)>0$. 

Note that $i(r) + j(r) = \Phi(r)$ where
$\Phi$ is the Hilbert polynomial of $\F$.

\medskip

{\bf Step 5: Structure of the hypothetical representing scheme $G^{(r)}$}
Let $\phi : \spec (R/I) \hra \spec R$ denote the inclusion and 
$\F'$ denote the pull-back of $\F$ under $\PP^n_{R/I} \hra \PP^n_R$.
The sheaf $\F'$ is flat over $R/I$, and the functor 
$\G^{(r)}_{R/I}$, which is the restriction of $\G^{(r)}$,
is represented by the linear scheme 
$\AA^d_{R/I} = \spec (R/I)[y_1,\ldots, y_d]$ over $R/I$,
where $d = \Phi(r)$ where $\Phi$ is the Hilbert polynomial of $\F$.
Hence, the pull-back of the hypothetical representing scheme $G^{(r)}$ to $R/I$
is the linear scheme $\AA^d_{R/I}$. 
We now use the following fact (see Lemma 6 and Lemma 7
of [N1] for a proof).

\medskip

{\bf Lemma }{\it
Let $R$ be a ring and $I$ a nilpotent ideal ($I^n =0$ for some $n\ge 1$). 
Let $X$ be a scheme over $\spec R$, such that the closed subscheme
$Y = X\otimes_R(R/I)$ is isomorphic over $R/I$
to $\spec B$ where $B$ is a finite-type $R/I$-algebra.
Let $b_1,\ldots, b_d\in B$ be a set of 
algebra generators for $B$ over $R/I$.  
Then $X$ is isomorphic over $R$ with $\spec A$ where $A$ is a finite-type
$R$-algebra. Moreover, there exists a set of 
$R$-algebra generators $a_1,\ldots, a_d$ for $A$, 
such that each $a_i$ restricts modulo $I$ to $b_i \in B$ over $R/I$.
Let $R[x_1,\ldots, x_d]$ be a polynomial ring in $d$ variables over
$R$, and consider the surjective $R$-algebra homomorphism 
$R[x_1,\ldots, x_d] \to A$ defined by sending each $x_i$ to $a_i$,
and let $J$ be its kernel. Then $J \subset I R[x_1,\ldots, x_d]$.
}

\medskip
It follows from the above lemma 
that $G^{(r)}$ is affine of finite type over $R$,
and its co-ordinate ring $A$ as an $R$ algebra is of the form
$$A = R[a_1,\ldots, a_d] = R[x_1,\ldots, x_d]/J$$ 
where $a_i$ is the residue of $x_i$, and 
$a_1,\ldots,a_d$ restrict over $R/I$ to the 
linear coordinates $y_1,\ldots, y_d$ 
on the linear scheme $\AA^d_{R/I}$,
and $J$ is an ideal with $J \subset I\cdot R[x_1,\ldots, x_d]$.
Being an additive group-scheme, $G^{(r)}$ has its zero section
$\sigma : \spec R \to G^{(r)}$, and this corresponds to 
an $R$-algebra homomorphism $\sigma^* : A \to R$.
Modulo $I$, the section $\sigma$ restricts to the 
zero section of $\AA^d_{R/I}$ over $\spec(R/I)$, therefore
$\sigma^*(a_i) \in I$ for all $i=1,\ldots,d$.  
Let $x'_i = x_i - \sigma^*(a_i) \in R[x_1,\ldots, x_d]$ and
$a'_i = a_i - \sigma^*(a_i) \in A$ be its residue modulo $J$.
Then $R[x_1,\ldots, x_d] = R[x'_1,\ldots, x'_d]$, the elements 
$a'_1, \ldots, a'_d$ generate $A$ as an $R$-algebra, and 
moreover the $a'_i$ restrict over $R/I$ to the 
linear coordinates $y_i$ 
on the linear scheme $\AA^d_{R/I}$. Therefore, by replacing the 
$x_i$ by the $x'_i$ and the
$a_i$ by the $a'_i$, we can assume that for each $i$, we have
$$\sigma^*(a_i) =0$$

Next, consider any element $f(x_1,\ldots, x_d) \in J$. Then
$f(a_1,\ldots a_d) = 0$ in $A$, so 
$\sigma^*f(a_1,\ldots a_d) = 0 \in R$, which shows that
the constant coefficient of $f$ is zero, as $\sigma^*(a_i) =0$. 
As we already know that $J \subset I\cdot R[x_1,\ldots, x_d]$, 
the vanishing of the constant term of any element of $J$ now
establishes that 
$$J \subset I\cdot (x_1,\ldots,x_d)$$
From the above, using $I^2=0$, 
it follows that for any $(b_1,\ldots,b_d) \in I^d$, we have a 
well-defined $R$-algebra homomorphism 
$$\Psi_{(b_1,\ldots,b_d)} : A \to R : a_i \mapsto b_i$$

We now express the linear-scheme structure of $G^{(r)}$ in 
terms of the ring $A$, using the fact that each $a_i$ restricts to 
$y_i$ modulo $I$, and $G^{(r)}_{R/I}$ is the standard linear-scheme 
$\AA^d_{R/I}$ with linear co-ordinates $y_i$. 
Note that the vector addition morphism
$\AA^d_{R/I}\times_{R/I}\AA^d_{R/I}\to \AA^d_{R/I}$
corresponds to the $R/I$-algebra homomorphism 
$$(R/I)[y_1,\ldots, y_d] \to 
(R/I)[y_1,\ldots, y_d] \otimes_{R/I} 
(R/I)[y_1,\ldots, y_d] : y_i \mapsto y_i\otimes 1 + 1\otimes y_i$$
while the scalar-multiplication morphism
$\AA^1_{R/I} \times_{R/I}  \AA^d_{R/I}\to \AA^d_{R/I}$
corresponds to the 
$R/I$-algebra homomorphism 
$$(R/I)[y_1,\ldots, y_d] \to (R/I)[t, y_1,\ldots, y_d] =
(R/I) [t]\otimes_{R/I} (R/I)[y_1,\ldots, y_d] :  
y_i \mapsto  ty_i$$

It follows that the addition morphism
$\alpha : G^{(r)} \times_R G^{(r)} \to G^{(r)}$ corresponds to
an algebra homomorphism $\alpha^* : A \to A\otimes_RA$ which has the form
$$a_i \mapsto a_i\otimes 1 + 1\otimes a_i + u_i 
\mbox{ where }u_i \in I(A\otimes_RA).$$
Let the element $u_i$ in the above equation for
$\alpha^*(a_i)$ be written as a polynomial expression
$$u_i = f_i(a_1\otimes 1, \ldots, a_d\otimes 1, 1\otimes a_1,
\ldots, 1\otimes a_d)$$ 
with coefficients in $I$. 
The additive identity $0$ of $G^{(r)}(R)$ corresponds to 
$\sigma^* : A \to R$
with $\sigma^*(a_i) =0$, and we have $0+0 =0$ in $G^{(r)}(R)$.
This implies that 
$f_i(0,\ldots,0) =0$, and so the constant term of $f_i$ is zero.
From this, using $I^2=0$, we get the important consequence that
$$f_i(w_1,\ldots,w_{2d}) = 0 \mbox{ for all }w_1,\ldots,w_{2d} \in I$$

The scalar-multiplication morphism
$\mu : \AA^1_R \times_R G^{(r)} \to G^{(r)}$ prolongs the 
standard scalar multiplication on $\AA^d_{R/I}$, and so 
$\mu$ corresponds to
an algebra homomorphism $\mu^* : A \to A[t] = R[t] \otimes_R A$ 
which has the form
$$a_i \mapsto ta_i + v_i \mbox{ where }v_i \in IA[t].$$
Let $v_i$ be expressed as a polynomial 
$v_i = g_i(t,a_1,\ldots,a_d)$ with coefficients in $I$. 
As multiplication by the scalar $0$ is the zero morphism on $G^{(r)}$,
it follows by specialising under $t\mapsto 0$ that 
$g_i(0, a_1,\ldots,a_d) =0$. This means $v_i = g_i(t,a_1,\ldots,a_d)$
can be expanded as a finite sum
$$v_i = \sum_{j \ge 1} t^j h_{i,j}(a_1,\ldots,a_d)$$
where the $h_{i,j}(a_1,\ldots,a_d)$ 
are polynomial expressions with coefficients 
in $I$. As the zero vector times any scalar is zero, 
it follows by specialising under $\sigma^*$ that 
$g_i(t, 0,\ldots,0) =0$. It follows that 
the constant term of each $h_{i,j}$ is zero. From this, and the fact that 
$I^2=0$, we get the important
consequence that 
$$g_i(t,b_1,\ldots,b_d) = 0 \mbox{ for all } b_1,\ldots,b_d \in I$$

{\bf Step 6: The kernel of the map $G^{(r)}(R) \to G^{(r)}(R/I)$}

{\bf Lemma }{\it 
Let $\Psi_{(b_1,\ldots,b_d)} : A \to R$
be the $R$-algebra homomorphism defined in terms of the generators by
$\Psi_{(b_1,\ldots,b_d)}(a_k) = b_k$. 
Let $\Psi : I^d \to Hom_{R-alg}(A,R)$ be the set-map defined by
$(b_1,\ldots,b_d) \mapsto (\Psi_{(b_1,\ldots,b_d)}:A\to R)$.
Then $\Psi$ is a homomorphism of $R$-modules, where 
the $R$-module structure on $Hom_{R-alg}(A,R)$ is defined by
its identification with the $R$-module $G^{(r)}(R)$.

The map $\Psi$ is injective, and its image is
the $R$-submodule $\ker G^{(r)}(\phi) \subset G^{(r)}(R)$, where 
$\phi : \spec (R/I) \hra \spec R$ is the inclusion.
}

\proof For any $(b_1,\ldots,b_d)$ and 
$(c_1,\ldots,c_d)$ in $I^d$, we have
\begin{eqnarray*}
& & (\Psi_{(b_1,\ldots,b_d)}+\Psi_{(c_1,\ldots,c_d)})(a_i)\\ 
& = & (\Psi_{(b_1,\ldots,b_d)}\otimes \Psi_{(c_1,\ldots,c_d)})
      (\alpha^*(a_i))\\
& = & b_i + c_i + f_i(b_1,\ldots, b_d, c_1,\ldots,c_d) 
\mbox { by substituting for }\alpha^*(a_i)\\
& = & b_i + c_i \mbox{ as } b_k, c_k \in I\\
& = & \Psi_{(b_1+c_1,\ldots,b_d+c_d)}(a_i).
\end{eqnarray*}
This shows the equality
$\Psi_{(b_1,\ldots,b_d)}+\Psi_{(c_1,\ldots,c_d)} =
\Psi_{(b_1,\ldots,b_d)+(c_1,\ldots,c_d)}$, 
which means the map $\Psi : I^d \to G^{(r)}(R)$ is additive.

For any $\lambda \in R$, let $f_{\lambda} : R[t] \to R$ be the 
$R$-algebra homomorphism defined by $f_{\lambda}(t) = \lambda$.
Then for any $(b_1,\ldots,b_d)\in I^d$ we have
\begin{eqnarray*}
(\lambda \cdot \Psi_{(b_1,\ldots,b_d)})(a_i) 
& = & (f_{\lambda} \otimes  \Psi_{(b_1,\ldots,b_d)})(\mu^*(a_i))\\
& = & (f_{\lambda} \otimes  \Psi_{(b_1,\ldots,b_d)})
      (ta_i + g_i(t,a_1,\ldots,a_d))\\
& = & \lambda b_i + g_i(\lambda, b_1,\ldots,b_d)\\
& = & \lambda b_i  \mbox{ as } b_k \in I\\
& = & \Psi_{(\lambda b_1,\ldots,\lambda b_d)}(a_i).
\end{eqnarray*}
This shows the equality
$\lambda \cdot \Psi_{(b_1,\ldots,b_d)}
= \Psi_{\lambda \cdot (b_1,\ldots,b_d)}$, hence 
the map $\Psi : I^d \to G^{(r)}(R)$ preserves scalar 
multiplication. This completes the proof that 
$\Psi : I^d \to G^{(r)}(R)$ is a homomorphism of $R$-modules.

The map $\Psi$ is clearly injective. The map
$G^{(r)}(\phi) : G^{(r)}(R)\to G^{(r)}(R/I)$ is
in algebraic terms the map 
$Hom_{R-alg}(A,R) \to Hom_{R-alg} (A,R/I)$
induced by the quotient $R \to R/I$. 
An element $g \in Hom_{R-alg} (A,R/I)$ represents the 
zero element of $G^{(r)}(R/I)$ exactly when
$g(a_i) =0 \in R/I$ for the generators $a_i$ of $A$.
Therefore $f \in Hom_{R-alg}(A,R)$ is in the kernel of $G^{(r)}(\phi)$ 
precisely when $f(a_i)\in I$ for the generators $a_i$.
Putting $b_i = f(a_i)$, we see that such an $f$ is
the same as $\Psi_{(b_1,\ldots,b_d)}$.

This completes the proof of the Lemma that $\ker G^{(r)}(\phi)= I^d$.

\medskip

In particular, as $\mm I=0$, it follows from the above Lemma 
that $\ker G^{(r)}(\phi)$ is annihilated by $\mm$,
so it is a vector space over $R/\mm$, and
its dimension as a vector space over $R/\mm$ is $d = \Phi(r)$, as by
assumption $I$ is a non-zero principal ideal.

The above determination of the 
dimension over $R/\mm$ of the kernel of $G^{(r)}(\phi)$ 
will contradict a more direct functorial description,
which is as follows.

\medskip

{\bf Step 7: Functorial description of kernel of $\G^{(r)}(R) \to
\G^{(r)}(R/I)$ }
As $\F_{R/I}(r)$ is flat over $R/I$, and as
for $r\ge m$ all higher direct images of $\F(r)$ vanish, 
$\G^{(r)}(R/I)$ is isomorphic to the $R/I$-module 
$(R/I)^d$ where $d=\Phi(r)$. 
By Lemma \ref{base change without flatness},
there exists $m'' \ge m$ such that for $r\ge m''$ the 
inclusion $\phi : \spec (R/I) \hra \spec R$ induces an isomorphism
$\phi^* \pi_*\, \F(r) \to \pi'_* \, \F'(r)$
where $\pi' : \PP^n_{R/I} \to \spec (R/I)$ is the projection and 
$\F'$ is the pull-back of $\F$ under $\PP^n_{R/I} \hra \PP^n_R$.
Note that $\G^{(r)}(R) = R^{i(r)}\oplus (R/I)^{j(r)}$,
and so for $r\ge m''$ we get an induced decomposition
$$\G^{(r)}(R/I) = (R/I)^{i(r)}\oplus (R/I)^{j(r)}$$ 
such that  
the map $\G^{(r)}(\phi) : \G^{(r)}(R) \to \G^{(r)}(R/I)$ is the 
map 
$$(q,1) : R^{i(r)}\oplus (R/I)^{j(r)} 
\to (R/I)^{i(r)}\oplus (R/I)^{j(r)}$$ 
where $q$ is the 
quotient map modulo $I$. It follows that the 
kernel of 
$\G^{(r)}(\phi)$ is the $R$-module 
$I^{i(r)}\oplus 0 \subset R^{i(r)}\oplus (R/I)^{j(r)} = \G^{(r)}(R)$. 
This is a vector
space over $R/\mm$ of dimension $i(r) < i(r)+j(r) = \Phi(r)$.

We thus obtain two different values for the dimension of 
the same vector space $\ker G^{(r)}(\phi) = \ker \G^{(r)}(\phi)$, 
which shows that our assumption that 
$\G^{(r)}$ is representable for arbitrarily large values of $r$ is false.
This completes the proof of the Theorem \ref{main theorem}.

\bigskip

{\bf Acknowledgement } This note was inspired by 
a question posed by Angelo Vistoli to the participants 
of the workshop `Advanced Basic Algebraic Geometry' held 
at the Abdus Salam ICTP, Trieste, in July 2003. 
The Corollary \ref{Question of Vistoli} answers that question.
I thank the ICTP for hospitality while this work was in progress.

\newpage

{\bf\large References}

\medskip

[EGA] Grothendieck, A. : {\it \'El\'ements de G\'eom\'etrie Alg\'ebriques}
(written with the collaboration of Dieudonn\'e). 
Publ. Math. IHES, {\bf 4}, {\bf 8}, {\bf 11}, {\bf 17}, 
{\bf 20}, {\bf 24}, {\bf 28}, {\bf 32} (1960-67).

[TDTE IV]  Grothendieck, A. : 
Techniques de construction et th\'eor\`emes d'existence
en g\'eom\'etrie alg\'ebriques IV : les sch\'emas de Hilbert.
S\'eminaire Bourbaki {\bf 221}, 1960/61.

[M] Mumford : {\it Lectures on Curves on an Algebraic Surface}. Princeton
University Press, 1966.

[N1] Nitsure, N. : Representability of ${\rm GL}_E$. 
Proc. Indian Acad. Sci. Math. Sci. {\bf 112} (2002). 
{\tt http://arXiv.org/abs/math/0204047}

[N2] Nitsure, N. : Construction of Hilbert and Quot Schemes.
Lectures at the ICTP Summer School on Advanced Basic Algebraic Geometry, 
July 2003. (Follow the web-link 
{\tt http://www.ictp.trieste/} 
for an electronic version).

\bigskip

\bigskip

{\small

\hfill School of Mathematics, 

\hfill Tata Institute of Fundamental Research,

\hfill Homi Bhabha Road, 

\hfill Mumbai 400 005, 

\hfill India. 

\hfill e-mail: {\tt nitsure@math.tifr.res.in}

}

\bigskip

\centerline{05 August 2003}

\end{document}